# Distributed Fault Detection and Accommodation in Dynamic Average Consensus

J. George, M. L. Elwin, R. A. Freeman, and K. M. Lynch

*Abstract*— This paper presents the formulation of fault detection and accommodation schemes for a network of autonomous agents running internal model-based dynamic average consensus algorithms. We focus on two types of consensus algorithms, one that is internally stable but non-robust to initial conditions and one that is robust to initial conditions but not internally stable. For each consensus algorithm, a fault detection filter based on the unknown input observer scheme is developed for precisely estimating the communication faults that occur on the network edges. We then propose a fault remediation scheme so that the agents could reach average consensus even in the presence of communication faults. Numerical results are provided to illustrate the efficacy of the proposed approach.

## I. INTRODUCTION

In networked systems, the individual agents are often required to reach consensus on certain quantities of interest. The problem of reaching consensus on the average of a set of local time-varying reference signals in a distributed fashion is typically known as *dynamic average consensus* and it plays a crucial role in several network applications such as distributed sensor fusion [1]–[3], formation control [4]–[6], and distributed mapping [7]–[9].

*Related Work:* Networks are often disrupted by faulty/malicious nodes, external attacks, or erroneous communication links. Although there exist consensus algorithms that can accommodate these network disruptions, almost all these approaches focus on static average consensus problems involving malicious nodes. These fault-or attack-tolerant approaches are called resilient consensus algorithms and they are mainly based on either a traditional fault accommodation scheme or simply disregard extreme values observed by the nodes. The latter method, commonly known as the mean-subsequence-reduced (MSR) algorithms, requires the nodes to reject the largest and smallest values it received from its neighbors and update their state based on the remaining values. An MSR-type algorithm was first proposed in [10] for a fully connected network and later extended to more general network topologies [11]. In [12], authors proposed a continuous-time version of the MSR algorithm to solve the asymptotic consensus problem when the maximum number of malicious nodes is known. The results of [12] were extended to dynamic networks in [13]. In [14] authors present an MSR-based resilient consensus algorithm for second-order, discrete-time systems. An extension of [14] for a known number of maximum malicious agents in the network is addressed in [15]. More recent developments on MSR-type algorithms can be found in [16] and [17]. Even though MSR-based algorithms have shown to be resilient to malicious nodes, they do not solve the fault detection problem and they are not concerned with erroneous communication scenarios. Also, currently there exists no straightforward way of extending these algorithms to a network of agents running dynamic average consensus protocols.

J. George is with the U.S. Army Research Laboratory, Adelphi, MD 20783, USA. jemin.george.civ@mail.mil

R. A. Freeman is with the Department of Electrical Engineering and Computer Science, Northwestern University, Evanston, IL 60208, USA. freeman@eecs.northwestern.edu

M. L. Elwin and K. M. Lynch are with the Department of Mechanical Engineering, Northwestern University, Evanston, IL 60208, USA. {elwin,kmlynch}@northwestern.edu

The problem of detecting and identifying misbehaving agents in a linear consensus network is first introduced in [18]. The unknown input estimator given in [18] assumes a single faulty node and requires knowledge of the global network structure. Several extensions of this unknown input observer based algorithm for the case of both Byzantine as well as non-colluding agents can be found in [19]–[21]. Further development of this algorithm for cyber-physical systems under attacks modeled by linear time-invariant descriptor systems with exogenous input is given in [22]. These works reveal that the necessary and sufficient conditions for the detectability of faulty nodes in a linear consensus network can be developed based on the strong detectability condition needed for the existence of the unknown input observers. In [23], authors propose a distributed algorithm based on the unknown input observers for the detection and isolation of faulty nodes. An extension of [23] to jointly detect and isolate faults occurring in nodes and edges is given in [24]. In [25], authors provide an explicit strategy that $k$ malicious nodes can follow to prevent a $2k$-connected node from computing the desired function of the initial state, or reaching consensus. Also shown in [25] is that the topology of the network completely characterizes the resiliency of the linear iterative algorithms and if the connectivity is $2k + 1$ or more, then the $k$ malicious nodes can be identified independent of their behavior. Almost all the unknown input observer based approaches require at least some nonlocal network structural information or measurements from non-neighbor nodes.

Besides the above-mentioned two main approaches, a distributed fault diagnosis architecture for large-scale dynamical systems is presented in [26]. In [27], a bank of distributed Kalman filters running a dynamic average consensus algorithm is used to tackle the problem of distributed fault diagnosis. A sliding mode observer based fault detection algorithm is given in [28] for the distributed detection of corrupted measurement exchange in a network of linear time-invariant systems. Though the algorithm in [28] successfully detects the presence of corrupted data using only local information, the fault isolation requires a centralized monitoring coordination layer. A clustering procedure based on both the similarity of measurements and the communication connectivity is proposed in [29] to address the problem of distributed fault detection and isolation of faulty nodes in the average consensus network. Since this clustering-based approach requires sequential processing of information, it scales poorly with network size.

*Contributions and Outline*

Almost all existing dynamic average consensus algorithms assume error-free communication and currently there exist no fault detection and accommodation schemes for dynamic average consensus algorithms with communication faults. Besides, the available distributed fault detection schemes rely on at least some global information. Therefore, we present the formulation of fault detection and accommodation schemes for a network running the internal model based dynamic average consensus algorithms given in [30]. Here we focus on two types of dynamic average consensus algorithms, one that is internally stable but non-robust to initial conditions and one that is robust to initial conditions but not internally stable. We develop an unknown input observer based fault detection filter for each algorithm and propose an approach that would allow the nodes to reach average consensus even in the presence of communication faults.

The structure of this paper is as follows. Details of the dynamic average consensus algorithms are first presented in section II. Formulation of the fault detection and accommodation scheme for the non-robust dynamic average consensus estimator and its numerical evaluation are given in sections III and IV, respectively. A fault detection and accommodation scheme for the robust dynamic average consensus estimator is given in section V. In section VI, we compare the performance of the proposed fault detection and accommodation schemes for the robust as well as non-robust estimator.

*Network Model*

Consider a connected undirected graph $\mathcal{G}(\mathcal{V}, \mathcal{E})$ of order $n$, where $\mathcal{V} \triangleq \{v_1, \ldots, v_n\}$ represents the nodes and $\mathcal{E} \subseteq \mathcal{V} \times \mathcal{V}$ represents the communication links between the nodes. The set of neighbors of node $v_i$ is denoted by $\mathcal{N}_i = \{v_j \in \mathcal{V} : (v_i, v_j) \in \mathcal{E}\}$. Let $\mathcal{A}$ be the *adjacency matrix* and $\Delta$ be the degree matrix associated with the graph $\mathcal{G}(\mathcal{V}, \mathcal{E})$. The graph *Laplacian* $\mathcal{L}$ induced by the information flow in $\mathcal{G}(\mathcal{V}, \mathcal{E})$ is defined as $\mathcal{L} = \Delta - \mathcal{A}$. The Laplacian's eigenvalues are listed in ascending order as $\lambda_1 \ldots \lambda_n$. The row sums of the Laplacian are always zero so $\lambda_1 = 0$ with corresponding eigenvector the vector of $n$ ones, $\mathbf{1}_n$. We denote the numerator and denominator polynomials of the transfer function $q(s)$ as $n_q(s)$ and $d_q(s)$, respectively.

## II. DYNAMIC AVERAGE CONSENSUS

This section presents two dynamic average consensus schemes based on the robust dynamic average consensus algorithm given in [30]. In the dynamic average consensus problem, all the nodes in the network have access to their own local reference signal and the nodes must estimate the average of all such signals within the network.

Let $\phi_i(t) \in \mathbb{R}$ denote the $i^{th}$ node's reference signal (input) at time $t$. Every node must track the time-varying signal:

$$\bar{\phi}(t) = \frac{1}{n}\mathbf{1}_n^T\mathbf{1}_n\,\phi(t),$$

where $\phi(t) = \begin{bmatrix}\phi_1(t) & \ldots & \phi_n(t)\end{bmatrix}^T$. We assume that each node's input, $\phi_i(s)$, has a transfer function of the following form:

$$\phi_i(s) = \frac{c_i(s)}{d(s)q_i(s)}, \quad (1)$$

where $c_i(s)$ is a polynomial that may differ between agents, $d(s)$ is a monic polynomial common to all agents with no roots in the open left half-plane and $q_i(s)$ is a monic polynomial with all of its roots in the open left half-plane. We assume that $c_i(s)$ and $d(s)q_i(s)$ are coprime. Examples of such $\phi_i(s)$ include the outputs of stable systems driven by step, ramp, sinusoidal, or exponentially growing signals. We further assume, without loss of generality, that $q_i(s) = 1$, because any open left half-plane pole of $\phi_i(s)$ contributes an exponentially vanishing component to $\phi_i(s)$.

We consider two dynamic average consensus estimators, related to the internal model (IM) estimator of [30], whose global dynamics are depicted in Figs. 1 and 2. We call the estimator in Fig. 1 the internally stable average consensus (ISAC) estimator and we call the estimator in Fig. 2 the robust average consensus (RAC) estimator. In both estimators, $h(s)$ and $g(s)$ correspond to the transfer function associated with the internal dynamics of the estimators and $\boldsymbol{\nu}(s) = \begin{bmatrix}\nu_1(s) & \cdots & \nu_n(s)\end{bmatrix}^T$ denotes the network estimate of $\bar{\phi}(t)$. Vectors $\boldsymbol{\eta}(s)$ and $\boldsymbol{\mu}(s)$ are auxiliary signals, where $\boldsymbol{\mu}(s)$ is kept internally for the ISAC estimator while agents communicate $\boldsymbol{\eta}(s)$ to their neighbors in the RAC estimator. Note that the implementation of these estimators only requires the communication of a single variable, i.e., $\boldsymbol{\nu}$ for the ISAC estimator and $\boldsymbol{\eta}$ for the RAC estimator. The effect of the initial conditions of the states of $h(s)$ and $g(s)$ appear as exogenous disturbances $\xi_h$ and $\xi_g$ in the block diagram. These disturbances have the same poles as their corresponding subsystems.

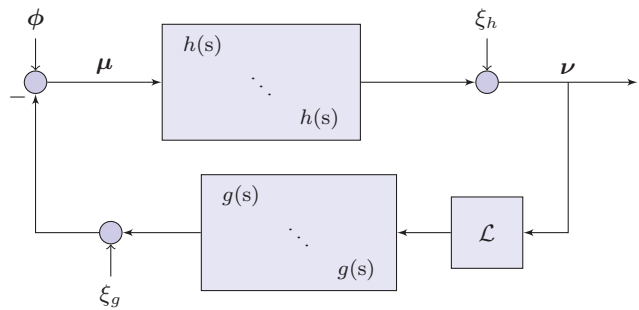

Fig. 1. Block diagram of the internally stable average consensus (ISAC) estimator; $\xi_h$ and $\xi_g$ represent exogenous inputs due to the initial conditions of the internal states of $h(s)$ and $g(s)$, respectively.

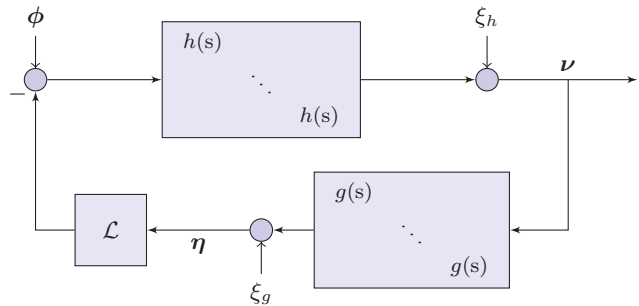

Fig. 2. Block diagram of the robust average consensus (RAC) estimator; $\xi_h$ and $\xi_g$ represent the exogenous inputs due to the initial conditions of the internal states of $h(s)$ and $g(s)$, respectively.

Both the ISAC and RAC estimators have the same transfer function from the input $\boldsymbol{\phi}$ to the average output $\boldsymbol{\nu}$:

$$\boldsymbol{\nu}(s) = (I + h(s)g(s)\mathcal{L})^{-1}h(s)\boldsymbol{\phi}(s). \quad (2)$$

However, their internal structure and their response to the initial condition disturbances differ; an important distinction when implementing the system on a network of nodes.

We first examine both estimators, assuming that the disturbances are zero (equivalently, the states of the $h(s)$ and $g(s)$ systems have been initialized to zero). The steady-state error for these estimators is given by

$$\mathbf{e}_{\text{ss}} = \lim_{t \to \infty}\left(\boldsymbol{\nu}(t) - \frac{1}{n}\mathbf{1}_n\mathbf{1}_n^T\boldsymbol{\phi}(t)\right). \quad (3)$$

The following theorem provides the conditions under which both of these estimators successfully track the average of the group's inputs with zero steady-state error.

*Lemma 1:* For a connected undirected graph $\mathcal{G}(\mathcal{V}, \mathcal{E})$, with $\xi_h = \xi_g = 0$ and $g(s)$ and $h(s)$ having no common unstable roots, the

ISAC and RAC estimators achieve dynamic average consensus with zero steady-state error (i.e., $\mathbf{e}_{ss} = \mathbf{0}$) if

i. ) the individual node inputs $\phi_i(s)$ satisfy (1) for all $i \in \{1, \ldots, n\}$,
ii. ) for some polynomial $p(s)$, $n_h(s) - d_h(s) = p(s)d(s)$ and $h(s)$ is stable,
iii. ) the roots of $d_g(s)d_h(s) + n_g(s)n_h(s)\lambda_i(\mathcal{L}) = 0$ are in the open left half-plane for all $i \in \{2, 3, \ldots, n\}$,
iv. ) for some polynomial $p_g(s)$, $d_g(s) = p_g(s)d(s)$.

*Proof:* For brevity, we provide an intuitive explanation, the details are similar to the proof in [30]. Conditions (i) and (iv) ensure zero steady-state error and arise from the internal model principle. Condition (iii) ensures that output components orthogonal to the consensus direction decay to zero (ensuring consensus). Condition (ii) ensures that when the agents reach consensus (which effectively disconnects the feedback path due to the Laplacian having an eigenvalue of 0 corresponding to the consensus direction $\mathbf{1}_n$) that the agents agree on the average. ∎

For the rest of this paper, we assume that $h(s)$ and $g(s)$ satisfy the restrictions set forth in Lemma 1.

Next, we examine the robustness of the estimators to initialization error. Robust estimators converge to the correct average regardless of the initial states of the $h(s)$ and $g(s)$ systems. Robustness provides resilience to situations when not all nodes start their estimators simultaneously, such as when nodes are added to or removed from the network or when the network splits into two separate connected components.

The following lemmas examine the robustness of the ISAC and RAC estimators with respect to the initial conditions of the $h(s)$ and $g(s)$ systems. We provide intuitive explanations and refer the reader to [30] for similar rigorous proofs.

*Lemma 2:* The ISAC and RAC estimators achieve zero steady-state tracking error ($\mathbf{e}_{ss} = \mathbf{0}$) for any initial conditions of the $h(s)$ system (i.e., $\xi_h \neq 0$).

*Proof:* This result follows because of the stability conditions (ii) and (iii), and because $\xi_h$ has the same poles as $h(s)$. ∎

*Lemma 3:* The ISAC estimator does not achieve zero steady-state tracking error for arbitrary initial conditions of $g(s)$.

*Proof:* This result occurs because, to satisfy Lemma 1, $g(s)$ and therefore the disturbance $\xi_g$ are not stable. Additionally, $\xi_g$ passes through to $\nu$ without first going through the Laplacian. See [31]. ∎

*Lemma 4:* The RAC estimator achieves zero steady-state tracking error for any initial conditions of $g(s)$.

*Proof:* For this estimator, $\xi_g$ passes through the Laplacian prior to reaching $\nu$, causing the effect of $\xi_g$ to dissipate. See [31]. ∎

The net result of Lemmas 2, 3, and 4 is that the RAC estimator is robust to initialization error, but the ISAC estimator is not. Another property that differs between these estimators is internal stability. As the following lemmas show, the ISAC estimator is internally stable, but the RAC estimator has an internal state that grows faster than the input and is not internally stable. The concern for internal stability arises from $g(s)$, which, because it contains the internal model $d(s)$ in its denominator, is not stable.

*Lemma 5:* The RAC estimator is not internally stable.

*Proof:* Neither the transfer function $g(s)$ nor its input are stable and the output, $\boldsymbol{\nu}$, reaches the integrator $g(s)$ without ever passing through a Laplacian [31]. ∎

*Lemma 6:* The ISAC estimator is internally stable.

*Proof:* Although $g(s)$ is not stable, when the estimator reaches consensus, the vector $\nu$ is parallel to $\mathbf{1}_n$, hence it is the eigenvector of $\mathcal{L}$ corresponding to the zero eigenvalue. Thus the input to $g(s)$ is zero. ∎

Overall, choosing the internally stable but non-robust ISAC estimator or the robust but not internally stable RAC estimator depends on the application. Both estimators, however, allow the nodes to converge to the average of the time-varying signals with zero steady-state error while communicating only one variable per agent. The estimator of [30] is both robust and internally stable; however, it requires every agent to communicate two values instead of just one. In the subsequent sections, we present fault detection and recovery schemes for both the ISAC and RAC estimators. We leave fault recovery for the estimator of [30] for future work.

## III. DISTRIBUTED FAULT DETECTION AND ACCOMMODATION

We now present a distributed fault detection and accommodation scheme for agents running the ISAC estimator. Here we focus on communication faults that are associated with edges of the network. Let $\mathcal{M}_i \subseteq \mathcal{N}_i$ be the subset of agent $i$'s neighbors from which agent $i$ receives faulty data. Any value received by agent $i$ from agent $j$ for all $j \in \mathcal{M}_i \subseteq \mathcal{N}_i$ is corrupted by an error signal $f_j^i$.

*Remark 1:* We do not assume any particular functional form for $f_j^i$ nor do we presume a known upper bound on $f_j^i$. The only assumption on $f_j^i$ is that it is piecewise differentiable.

Implementation of the traditional internal model based dynamic average consensus estimator requires the communication of both $\nu$ and $\eta$ between two adjacent nodes. The ISAC and RAC estimators given here require the nodes to communicate only one output, i.e., $\nu$ for the ISAC estimator and $\eta$ for the RAC estimator. Furthermore, the traditional consensus estimator assumes perfect communication between nodes. In reality, communication links are imperfect and any error in communication would prevent the nodes from reaching consensus. Here we show that if we permit the nodes to communicate a second variable besides $\nu$, then the nodes would be able to detect any communication faults in its incoming links and deploy schemes to counteract the adverse effects of communication errors.

### A. Fault Detection Filter

To simplify notation and facilitate the formulation of the fault detection and isolation scheme, we assume that the transfer function $h(s)$ is strictly proper. Note that either $h(s)$ or $g(s)$ must be strictly proper to implement ISAC.

The state equations that every agent executes for the ISAC algorithm are

$$\dot{X}_i^1 = A_1 X_i^1 - B_1 \mu_i, \tag{4a}$$
$$\nu_i = C_1 X_i^1, \tag{4b}$$

and

$$\dot{X}_i^2 = A_2 X_i^2 + B_2 \sum_{j \in \mathcal{N}_i} (\nu_i - \nu_j), \tag{5a}$$
$$\mu_i = C_2 X_i^2 + D_2 \sum_{j \in \mathcal{N}_i} (\nu_i - \nu_j) - \phi_i. \tag{5b}$$

Here $X_i^1 \in \mathbb{R}^{m_1}$ is agent $i$'s state corresponding to the $h(s)$ subsystem and $X_i^2 \in \mathbb{R}^{m_2}$ is agent $i$'s state corresponding to the $g(s)$ subsystem, from the block diagram in Fig. 1. Agent $j$ communicates $\nu_j$ to its neighbors.

To facilitate fault detection and recovery, agents communicate $\mu_i$ in addition to $\nu_i$. Therefore, agent $i$ receives the following corrupt signals from agent $\ell$ for all $\ell \in \mathcal{M}_i \subseteq \mathcal{N}_i$:

$$\widetilde{\nu}_\ell^i = C_1 X_\ell^1 + f_\ell^i, \tag{6}$$

$$\widetilde{\mu}_\ell^i = C_2 X_\ell^2 + D_2 \sum_{j \in \mathcal{N}_\ell} (\nu_\ell - \nu_j) - \phi_\ell + f_\ell^i. \tag{7}$$

Consider the following extended system, representing agent $i$'s view of agent $\ell$'s $X^1$ dynamics:

$$\begin{bmatrix} \dot{X}_\ell^1 \\ \dot{f}_\ell^i \end{bmatrix} = \begin{bmatrix} A_1 & B_1 \\ \mathbf{0}_{m_1 \times 1}^T & 0 \end{bmatrix} \begin{bmatrix} X_\ell^1 \\ f_\ell^i \end{bmatrix} + \begin{bmatrix} -B_1 \\ 0 \end{bmatrix} \widetilde{\mu}_\ell^i + \begin{bmatrix} \mathbf{0}_{m_1 \times 1} \\ 1 \end{bmatrix} \dot{f}_\ell^i, \tag{8}$$

$$\widetilde{\nu}_\ell^i = \begin{bmatrix} C_1 & 1 \end{bmatrix} \begin{bmatrix} X_\ell^1 \\ f_\ell^i \end{bmatrix}. \tag{9}$$

The above extended system can be written as

$$\dot{\mathbf{x}}_\ell^i = A\mathbf{x}_\ell^i + B\widetilde{\mu}_\ell^i + E\dot{f}_\ell^i, \tag{10a}$$

$$\widetilde{\nu}_\ell^i = C\mathbf{x}_\ell^i, \tag{10b}$$

where $\mathbf{x}_\ell^i = \begin{bmatrix} X_\ell^1 \\ f_\ell^i \end{bmatrix}$, $A = \begin{bmatrix} A_1 & B_1 \\ 0 & 0 \end{bmatrix}$, $B = \begin{bmatrix} -B_1 \\ 0 \end{bmatrix}$, $E = \begin{bmatrix} 0 \\ 1 \end{bmatrix}$, and $C = \begin{bmatrix} C_1 & 1 \end{bmatrix}$. Based on the observations $\widetilde{\nu}_\ell^i$ and $\widetilde{\mu}_\ell^i$, node $i$ deploys the following full-order observer:

$$\dot{\mathbf{z}}_\ell^i = F\mathbf{z}_\ell^i + TB\widetilde{\mu}_\ell^i + K\widetilde{\nu}_\ell^i, \tag{11a}$$

$$\widehat{\mathbf{x}}_\ell^i = \mathbf{z}_\ell^i + H\widetilde{\nu}_\ell^i, \tag{11b}$$

where $F, K, T$, and $H$ are matrices to be designed for achieving desirable estimation error $(\widetilde{\mathbf{x}}_\ell^i = \mathbf{x}_\ell^i - \widehat{\mathbf{x}}_\ell^i)$ performance as shown next.

*Theorem 1:* For the extended system given in (10) and the observer given in (11), the estimation error $\widetilde{\mathbf{x}}_\ell^i$ is exponentially stable if the matrices $F$, $H$, $T$, and $K = K_1 + K_2$ are selected such that

$$(HC - I) E = 0, \tag{12}$$

$$T = I - HC, \tag{13}$$

$$F = (A - HCA - K_1 C), \tag{14}$$

$$K_2 = FH, \tag{15}$$

and $K_1$ is selected such that $F$ is strictly stable.

*Proof:* After some algebraic manipulations, the estimation error dynamics can be written as

$$\dot{\widetilde{\mathbf{x}}}_\ell^i = (A - HCA - K_1 C) \widetilde{\mathbf{x}}_\ell^i - [F - (A - HCA - K_1 C)] \mathbf{z}_\ell^i \\ - [K_2 - (A - HCA - K_1 C) H] \widetilde{\nu}_\ell^i \\ - [T - (I - HC)] B\widetilde{\mu}_\ell^i - (HC - I) E\dot{f}_\ell^i, \tag{16}$$

Substituting (12), (13), (14), and (15) yields

$$\dot{\widetilde{\mathbf{x}}}_\ell^i = F\widetilde{\mathbf{x}}_\ell^i. \tag{17}$$

If all eigenvalues of $F$ have negative real parts, then the estimator is exponentially stable. ∎

*Remark 2:* The observer in (11) is typically known as the unknown input observer because the estimation error asymptotically approaches zero, regardless of the presence of the unknown input $\dot{f}_\ell^i$. The necessary and sufficient conditions for the existence of such an observer are given next.

*Theorem 2:* For the extended system given in (10) and the observer given in (11), the estimation error $\widetilde{\mathbf{x}}_\ell^i$ is exponentially stable if and only if

(i) $\text{rank}(CE) = \text{rank}(E)$, and

(ii) invariant zeros of $(A, E, C)$ are strictly stable.

*Proof:* A proof and more details on unknown input observers can be found in [32] and [33]. ∎

If condition (i) is satisfied, then a matrix $H$ that satisfies (12) can be calculated as

$$H = E \left[ (CE)^T (CE) \right]^{-1} (CE)^T. \tag{18}$$

Furthermore, if $(A - HCA, C)$ is a detectable pair, then a gain matrix $K_1$ can be easily obtained such that $F$ is Hurwitz.

*Remark 3:* Therefore, after implementing an observer for each of its neighbors, node $i$ would be able to precisely estimate any faults in its incoming links. Note that the implementation of the proposed observer does not require any nonlocal information or knowledge of the network structure.

Given next is a fault recovery scheme that allows the nodes to accommodate any communication faults and reach consensus despite the presence of any such faults.

### B. Fault Accommodation

If, for each potentially faulty incoming link, the agents implement a fault detection filter from the previous section, then they can compensate for any faults using the scheme we develop here. For illustrative purposes, assume agent $i$ receives corrupt information from agent $\ell$ and all other links are error-free; however, our scheme works with multiple corrupt communication links. Let $\widehat{f}_\ell^i$ be the estimates of $f_\ell^i$ obtained by the fault detection filter given in (11). Since the estimation error $\widetilde{\mathbf{x}}_\ell^i$ is guaranteed to be exponentially stable, node $i$ can account for the erroneous $\widetilde{\nu}_\ell^i$ by subtracting the estimated faults from the erroneous signal. Thus the new dynamic average consensus algorithm can be written as

$$\dot{X}_i^1 = A_1 X_i^1 - B_1 \mu_i, \tag{19}$$

$$\dot{X}_i^2 = A_2 X_i^2 + B_2 \left( \sum_{j \in \mathcal{N}_i \setminus v_\ell} (\nu_i - \nu_j) + \left( \nu_i - \widetilde{\nu}_\ell^i + \widehat{f}_\ell^i \right) \right) \tag{20}$$

$$\nu_i = C_1 X_i^1, \tag{21}$$

$$\mu_i = C_2 X_i^2 + D_2 \left( \sum_{j \in \mathcal{N}_i \setminus v_\ell} (\nu_i - \nu_j) + \left( \nu_i - \widetilde{\nu}_\ell^i + \widehat{f}_\ell^i \right) \right) - \phi_i. \tag{22}$$

*Remark 4:* Although the fault estimate converges asymptotically to the actual fault, the transient error in this estimate enters the system in the same manner as an initial condition for the $g(s)$ subsystem. Therefore, unless the fault estimator state is initialized to zero, there will be a transient in the fault estimate and thus a non-zero steady-state error (see Lemma 3). Our fault accommodation scheme for the RAC estimator does not suffer from this drawback.

### IV. EXAMPLE I

We evaluate the performance of the ISAC fault detection and accommodation scheme through numerical simulation. Consider the connected undirected graph of 9 nodes shown in Fig. 3.

The inputs $\phi_i(t)$ are

$$\phi_i(t) = i \sin \left( \omega t + \frac{i}{4} \pi \right), \quad \text{for all } i \in \{1, 2, 3, 4, 5\}, \\ \phi_i(t) = i \cos \left( \omega t + \frac{i}{4} \pi \right), \quad \text{for all } i \in \{6, 7, 8, 9\}. \tag{23}$$

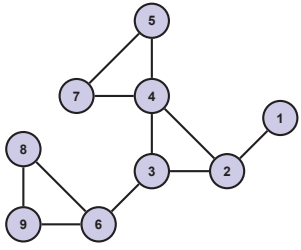

Fig. 3. Network topology.

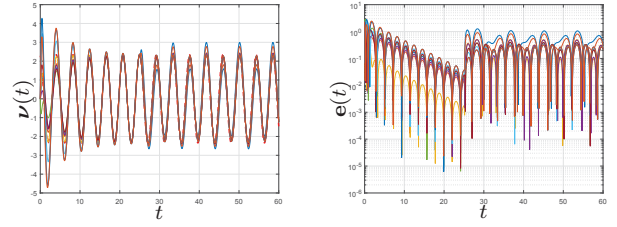

(a) Consensus estimates

(b) Consensus error

Fig. 5. Consensus estimates and errors for the erroneous scenario.

where the frequency $\omega = 1.5$. Here $d(s) = (s^2 + \omega^2)$ (sinusoids with frequency $\omega$). We choose $h(s)$ to be

$$h(s) = \frac{2\omega s + 3\omega^2}{s^2 + 2\omega s + 4\omega^2}$$

and $g(s)$ to be

$$g(s) = \frac{1.5s}{s^2 + \omega^2}.$$

It can be verified that $h(s)$ and $g(s)$ satisfy the conditions for Lemma 1 for our chosen graph. The matrices $A_1, A_2, B_1, B_2, C_1,$ and $C_2$ can now be obtained from $h(s)$ and $g(s)$. Note that $D_1 = D_2 = 0$ for the selected $h(s)$ and $g(s)$.

Figure 4 displays the results of running the ISAC estimator in the absence of communication faults.

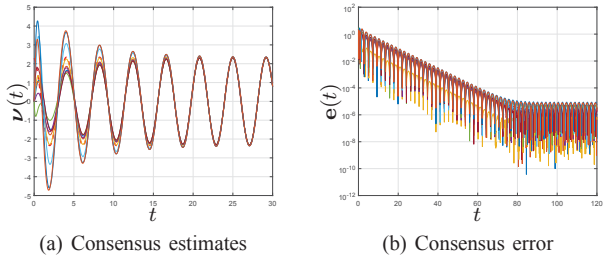

(a) Consensus estimates

(b) Consensus error

Fig. 4. ISAC estimates and estimation error for faultless communication.

Now we consider a scenario where the following communication error is associated with the link between agents 1 and 2.

$$f_1^2(t) = f_2^1(t) = \begin{cases} 0, & \text{if } t \leq 25; \\ \cos\left(\tfrac{1}{2}\omega t\right), & \text{otherwise.} \end{cases} \quad (24)$$

**Remark 5:** For ease of simulation, here we assume the fault is symmetric, i.e., $f_1^2(t) = f_2^1(t)$. The proposed fault detection and accommodation scheme works for non-symmetric faults. Also, the sinusoidal fault signal is selected because the individual reference signals are sinusoidal and thus it may be difficult to detect a sinusoidal fault signal rather than an arbitrary fault. As shown in Theorem 1, the proposed fault detection and accommodation scheme works for any type of faults as along as it is piecewise differentiable.

The results obtained by the ISAC estimator without fault accommodation are given in Fig. 5 and show that the fault introduces error into the estimate.

After constructing the extended system given in (10), we have $C = \begin{bmatrix} 3 & 6.75 & 1 \end{bmatrix}$ and $E = \begin{bmatrix} 0 & 0 & 1 \end{bmatrix}^T$. Note that $\text{rank}(CE) = \text{rank}(E) = 1$ and a matrix $H$ that satisfies (12) can be selected as $H = E$. The pair $(A - HCA, C)$ is detectable. The gain matrix $K_1$ is $\begin{bmatrix} 5.3993 & 12.1485 & 1.7998 \end{bmatrix}^T$. Now matrices $F$, $T$, and $K_2$ can be calculated from (13), (14), and (15), respectively.

The estimation errors $\widetilde{\mathbf{x}}_2^1$ and $\widetilde{\mathbf{x}}_1^2$, as well as the estimated faults $\hat{f}_2^1$ and $\hat{f}_1^2$ obtained from implementing the fault detection filter given in (11a) and (11b) for the erroneous scenario, are given in Fig. 6.

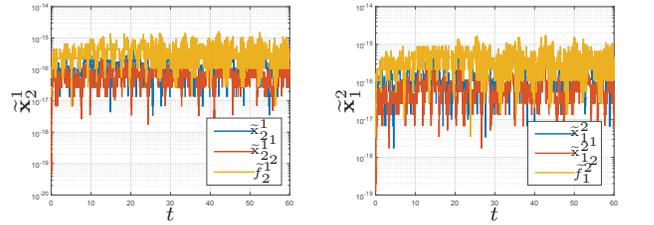

(a) Fault detection filter for node 1

(b) Fault detection filter for node 2

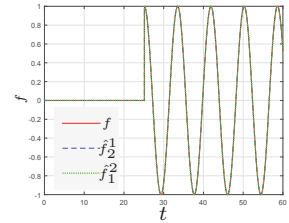

(c) Estimated communication faults

Fig. 6. Fault detection filter estimation errors and estimated communication faults

Finally, the dynamic average consensus algorithm given in subsection III-B is implemented using the error estimates obtained from the fault detection filters. Figure 7 shows that, when using our fault accommodation scheme, the system converges with zero steady-state error despite the presence of communication faults.

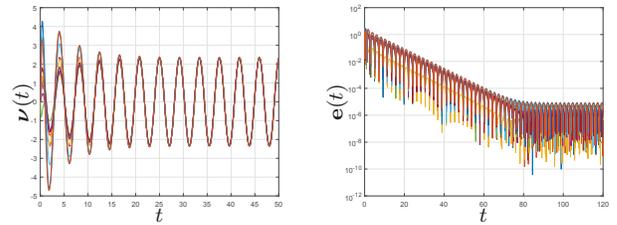

(a) Consensus estimates

(b) Consensus error

Fig. 7. Consensus estimates and corresponding errors for the ISAC with the fault detection and accommodation algorithm.

## V. ROBUST FAULT DETECTION AND ACCOMMODATION

We now present a distributed fault detection and accommodation scheme for agents running the RAC estimator. We use the same fault model as in section III, with faults being associated with communication links. Unlike the ISAC estimator and the associated fault compensation scheme, this estimator is robust to initial condition errors. This robustness, however, comes at the expense of internal stability.

### A. Robust Fault Detection Filter

The dynamics of the individual agents implementing the RAC estimator can be written as

$$\dot{X}_i^1 = A_1 X_i^1 + B_1 \left( \phi_i - \sum_{j \in \mathcal{N}_i} (\eta_i - \eta_j) \right), \quad (25a)$$

$$\nu_i = C_1 X_i^1 + D_1 \left( \phi_i - \sum_{j \in \mathcal{N}_i} (\eta_i - \eta_j) \right), \quad (25b)$$

and

$$\dot{X}_i^2 = A_2 X_i^2 + B_2 \nu_i, \quad (26a)$$
$$\eta_i = C_2 X_i^2 \quad (26b)$$

where $X_i^1 \in \mathbb{R}^{m_1}$ is the state associated with the $h(s)$ subsystem and $X_i^2 \in \mathbb{R}^{m_2}$ is the state associated with the $g(s)$ subsystem in the block diagram of Fig. 2. We have assumed that $g(s)$ is strictly proper, so $D_2 = 0$.

To implement RAC, the agents communicate the intermediate output $\eta_i$ rather than their estimate $\nu_i$. To implement fault detection, the agents also must communicate $\nu_i$.

When the communication channel from agent $\ell$ to agent $i$ is corrupted by $f_\ell^i$, for all $\ell \in \mathcal{M}_i \subseteq \mathcal{N}_i$, agent $i$ receives the following observations

$$\widetilde{\nu}_\ell^i = -\nu_\ell + f_\ell^i, \quad (27)$$
$$\widetilde{\eta}_\ell^i = \eta_\ell + f_\ell^i. \quad (28)$$

Now consider the following extended system, which is agent $i$'s view of agent $\ell$'s $X^2$ dynamics:

$$\begin{bmatrix} \dot{X}_\ell^2 \\ \dot{f}_\ell^i \end{bmatrix} = \begin{bmatrix} A_2 & B_2 \\ \mathbf{0}_{m_2 \times 1}^T & 0 \end{bmatrix} \begin{bmatrix} X_\ell^2 \\ f_\ell^i \end{bmatrix} + \begin{bmatrix} -B_2 \\ 0 \end{bmatrix} \widetilde{\nu}_\ell^i + \begin{bmatrix} \mathbf{0}_{m_2 \times 1} \\ 1 \end{bmatrix} \dot{f}_\ell^i, \quad (29)$$

$$\widetilde{\eta}_\ell^i = \begin{bmatrix} C_2 & 1 \end{bmatrix} \begin{bmatrix} X_\ell^2 \\ f_\ell^i \end{bmatrix}. \quad (30)$$

After defining $\mathbf{x}_\ell^i = \begin{bmatrix} X_\ell^2 \\ f_\ell^i \end{bmatrix}$, $A = \begin{bmatrix} A_2 & B_2 \\ 0 & 0 \end{bmatrix}$, $B = \begin{bmatrix} -B_2 \\ 0 \end{bmatrix}$, $E = \begin{bmatrix} 0 \\ 1 \end{bmatrix}$, and $C = \begin{bmatrix} C_2 & 1 \end{bmatrix}$, the above extended system can be placed into the same form as the system (10). Now based on the observations $\widetilde{\eta}_\ell^i$ and $\widetilde{\nu}_\ell^i$, node $i$ deploys the following full-order observer:

$$\dot{\mathbf{z}}_\ell^i = F \mathbf{z}_\ell^i + T B \widetilde{\nu}_\ell^i + K \widetilde{\eta}_\ell^i, \quad (31a)$$
$$\widehat{\mathbf{x}}_\ell^i = \mathbf{z}_\ell^i + H \widetilde{\eta}_\ell^i, \quad (31b)$$

where the matrices $F, K, T$, and $H$ are selected such that conditions (12), (13), (14), and (15) are satisfied. Now following the same reasoning given in Theorem 2, we can conclude that the observer given in (31) is exponentially stable.

### B. Robust Fault Accommodation

As with the ISAC fault accommodation scheme, the RAC scheme assumes that every agent implements a fault detection filter for each potentially faulty incoming communication link. For illustrative purposes, we assume agent $i$ receives corrupt information from agent $\ell$ and the other communication links have no faults; however, this scheme extends to multiple faulty links. Let $\widehat{f}_\ell^i$ be the estimate of $f_\ell^i$ obtained by the fault detection filter described in section V-A. Since the estimation error $\widetilde{\mathbf{x}}_\ell^i$ is guaranteed to be exponentially stable, node $i$ can account for the erroneous $\widetilde{\eta}_\ell^i$ by subtracting the estimated faults from the erroneous signal. Thus the new dynamic average consensus algorithm can be written as

$$\dot{X}_i^1 = A_1 X_i^1 + B_1 \left( \phi_i - \sum_{j \in \mathcal{N}_i \setminus v_\ell} (\eta_i - \eta_j) - \left( \eta_i - \widetilde{\eta}_\ell^i + \widehat{f}_\ell^i \right) \right) \quad (32a)$$

$$\nu_i = C_1 X_i^1 + D_1 \left( \phi_i - \sum_{j \in \mathcal{N}_i \setminus v_\ell} (\eta_i - \eta_j) - \left( \eta_i - \widetilde{\eta}_\ell^i + \widehat{f}_\ell^i \right) \right) \quad (32b)$$

and

$$\dot{X}_i^2 = A_2 X_i^2 + B_2 \nu_i, \quad (33a)$$
$$\eta_i = C_2 X_i^2. \quad (33b)$$

***Remark 6:*** Since the transient error in the fault estimate enters the system in the same manner as an initial condition for the $h(s)$ subsystem, the proposed fault accommodation scheme for the RAC estimator achieves zero steady-state tracking error ($\mathbf{e}_{ss} = \mathbf{0}$) for any initial conditions of the fault estimator state.

## VI. EXAMPLE II

Here we evaluate the performance of the fault detection and accommodation algorithm for the RAC estimator though numerical simulations. We consider the same connected undirected graph of 9 nodes given in Fig. 3. The inputs $\phi_i(t)$ are the same as those given in (23). Thus the same $h(s)$ and $g(s)$ from the previous example can be used to implement the robust dynamic average consensus algorithm.

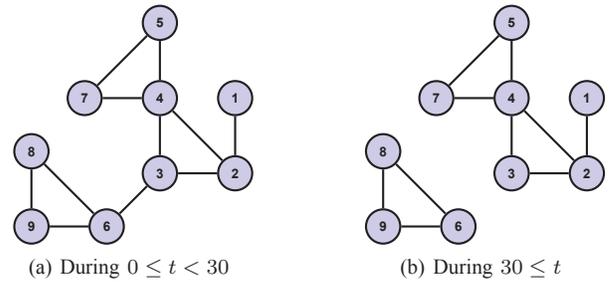

(a) During $0 \leq t < 30$   (b) During $30 \leq t$

Fig. 8. Network topology before (a) and after (b) the split.

We demonstrate the robustness of the dynamic average consensus algorithms by severing the link between agents 3 and 6 at time $t = 30$, which splits the graph into two disjoint connected components for the rest of the simulation time (see Fig. 8). Thus for $t \geq 30$, nodes $v_1, v_2, v_3, v_4, v_5, v_7$ should converge to

$$\bar{\phi}_1(t) = \frac{1}{6} \left( \sum_{i=1}^{5} \phi_i(t) + \phi_7(t) \right),$$

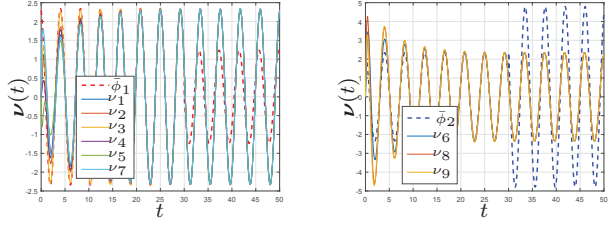

(a) ISAC estimates for nodes $v_1, v_2, v_3, v_4, v_5, v_7$   (b) ISAC estimates for nodes $v_6, v_8, v_9$

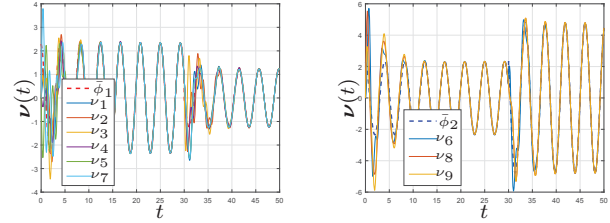

(c) RAC estimates for nodes $v_1, v_2, v_3, v_4, v_5, v_7$   (d) RAC estimates for nodes $v_6, v_8, v_9$

Fig. 9. Consensus estimates for the error-free scenario. The network splits at $t = 30$.

and nodes $v_6$, $v_8$, $v_9$ should converge to

$$\bar{\phi}_2(t) = \frac{1}{3}\left(\phi_6(t) + \phi_8(t) + \phi_9(t)\right).$$

The results obtained from both the RAC and ISAC estimators without communication error are shown in Figs. 9 and 10. Due to the non-robustness of the ISAC estimator, when the graph splits, the agents in each connected component converge to the incorrect average (see Figs. 9(a) and 9(b)). The robustness of the RAC estimator allows the agents in each connected component to reach the average of the inputs in their respective subgraphs. However, this robustness comes at the price of internal stability, so there is a tradeoff between these two estimators.

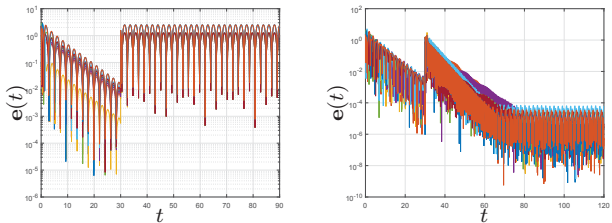

(a) ISAC estimator error   (b) RAC estimator error

Fig. 10. Consensus errors for the error-free scenario. The network splits at $t = 30$.

Next, we consider a scenario where the communication error between agents 1 and 2 is given in (24). When the agents implement the RAC estimator without fault accommodation, their estimates are erroneous, as shown in Fig. 11. Note that the fault first becomes nonzero at $t = 25$, which causes some initial error. When the network splits at $t = 30$, the component of the network that did not contain the faulty link converges to the desired value, while the component containing the faulty link continues to have an error.

We implement the fault detection filter as described in section V-A, based on (31). The estimated faults $f_2^1$ and $f_1^2$ from the fault detection filter for the erroneous scenario are shown in Fig. 12.

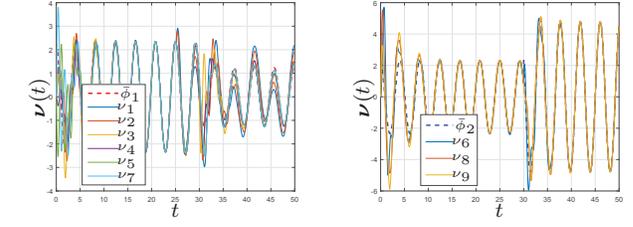

(a) RAC estimate for nodes $v_1, v_2, v_3, v_4, v_5, v_7$   (b) RAC estimates for nodes $v_6, v_8, v_9$

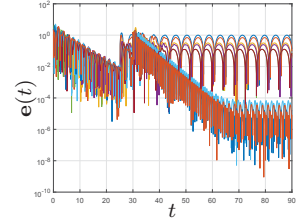

(c) RAC consensus error

Fig. 11. RAC estimates and error. The communication fault begins at $t = 25$ and the networks splits at $t = 30$.

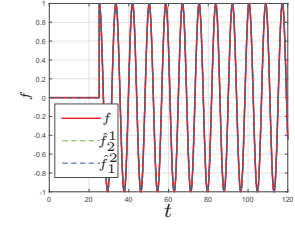

Fig. 12. Communication fault estimates

Finally, the dynamic average consensus algorithm given in subsection V-B is implemented using the error estimates obtained from the fault detection filters. Results obtained from implementing the robust fault accommodation scheme are given in Fig. 13. These results show that even with the fault and the network split, the agents' estimates converge to the correct average.

## VII. CONCLUSION

We have presented fault detection and recovery schemes for two types of dynamic average consensus estimators, one that is internally stable but non-robust to initial conditions (ISAC) and one that is robust to initial conditions but not internally stable (RAC). The agents use unknown input observers to estimate their neighbor's state and any faults on the communication links. These estimates are then used to recover from the fault and reach average consensus. Future work includes extending these methods to a robust and internally stable average consensus estimator. We also plan to develop stochastic version the proposed algorithm to account for random communication noise.

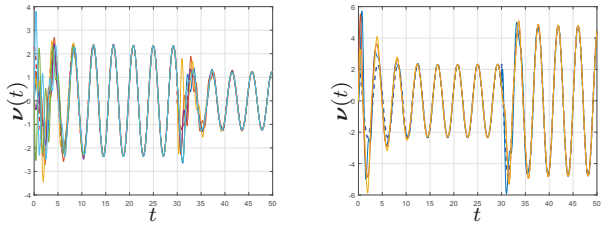

(a) RAC estimates with fault recovery for nodes $v_1, v_2, v_3, v_4, v_5, v_7$

(b) RAC estimates with fault recovery for nodes $v_6, v_8, v_9$

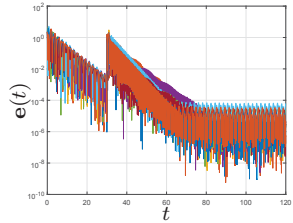

(c) Robust consensus error

Fig. 13. RAC estimates and corresponding errors from fault accommodating RAC estimator. The fault begins at $t = 25$ and the networks splits at $t = 30$.